\documentclass[reqno]{amsart}
\RequirePackage[english]{babel}
\RequirePackage{amssymb}
\usepackage{amsmath,amsfonts,amsthm}
\usepackage{mathtools}
\RequirePackage{xcolor}
\definecolor{linkblue}{rgb}{0 0.2 0.6}
\RequirePackage{hyperref}
\hypersetup{
	colorlinks=true,
	pdfborder={0 0 1},
	linkcolor=linkblue,
	citecolor=linkblue,
	urlcolor=linkblue,
	breaklinks,
}

\newtheorem{thm}{Theorem}
\newtheorem{lem}[thm]{Lemma}
\newtheorem{cor}[thm]{Corollary}

\theoremstyle{remark}

\theoremstyle{definition}
\newtheorem{dfn}[thm]{Definition}

\DeclareMathOperator{\M}{\mathsf E}
\DeclareMathOperator{\cov}{\mathbf{cov}}

\newcommand{\prob}{\stackrel{\text{\rm P}}{\rightarrow}}

\newcommand{\convdistr}{\stackrel{\text{\rm d}}{\rightarrow}}

\numberwithin{equation}{section}

\begin{document}
	\renewcommand{\theequation}{\thesection.\arabic{equation}}
	\title
	[Asymptotic normality of EW-TLS estimator in a multivariate EIV model]
	{Asymptotic normality of element-wise weighted total least squares estimator in a multivariate errors-in-variables model}
	
	\author{Ya. V. Tsaregorodtsev}
	\address{Department of Mathematical Analysis,
		Faculty of Mechanics and Mathematics,
		Taras Shevchenko National University of Kyiv,
		Building 4-e, Akademika Glushkova Avenue,
		Kyiv, Ukraine, 03127}
	\email{777Tsar777@mail.ru}
	
	
	\begin{abstract}
		A multivariable measurement error model $AX \approx B$ is considered. Here $A$ and $B$ are input and output matrices of measurements and $X$ is a rectangular matrix of fixed size to be estimated. The errors in $[A,B]$ are row-wise independent, but within each row the errors may be correlated. Some of the columns are observed without errors and the error covariance matrices may differ from row to row. The total covariance structure of the errors is known up to a scalar factor. The fully weighted total least squares estimator of $X$ is studied. We give conditions for asymptotic normality of the estimator, as the number of rows in $A$ is increasing. We provide that the covariance structure of the limiting Gaussian random matrix is nonsingular. 
	\end{abstract}
	
	\maketitle
	
	\section{Introduction}\label{s:1}
	We deal with an overdetermined set of linear equations $AX \approx B,$ which is common in linear parameter estimation problems \cite{huva}. If both the data matrix $A$ and observation matrix $B$ are contaminated with errors, and all the errors are uncorrelated and have equal variances, the total least squares (TLS) technique
	is appropriate for solving this set \cite{golo}, \cite{huva}. Under mild conditions, the TLS estimator of $X$ is consistent and asymptotically normal, as the number of rows in $A$ is increasing \cite{gl}, \cite{kutsa}.
	
	In this paper we consider heteroscedastic errors. The errors in $[A, B]$ are row-wise independent, but within each row the errors may be correlated. Some of the columns are observed without errors, and the error covariance matrices may differ from row to row. The total error covariance structure is assumed known up to a scalar factor. For this model, the element-wise weighted total least squares (EW-TLS) estimator is introduced and its consistency is proven in \cite{kuhu}. Concerning the computation of the estimator see \cite{mark}, \cite{jaam}. The EW-TLS estimator $\hat{X}$ is applied, e.g.,  in geodesy \cite{mah}. 
	
	Our goal is to extend the asymptotic normality result of \cite{kutsa} to the EW-TLS estimator. We work under the conditions of Theorem 2, \cite{kuhu} about the consistency of $\hat{X}.$ We use the objective function of the estimator, see formula (22) in \cite{kuhu}, and the rules of matrix calculus~\cite{cart}.
	
	The paper is organized as follows. In section \ref{s:2}, we describe the model, introduce main assumptions, refer to the consistency result for $\hat{X}$ and present the objective function and the matrix estimating function. In Section \ref{s:3}, we state the asymptotic normality result and provide a nonsingular covariance structure for a limiting random matrix. In Section \ref{s:4}, we derive consistent estimators for nuisance parameters of the model in order  to estimate consistently the asymptotic covariance structure of $\hat{X},$ and Section \ref{s:5} concludes. The proofs are given in  Appendix. 
	
	Throughout the paper all vectors are column ones, $\M$ stands for expectation and acts as an operator  on the total product, $\cov(x)$ denotes the covariance matrix of a random vector $x,$ and for a sequence of random matrices $\{X_m, m\geq 1\}$ of the same size, notation $X_m = O_p(1)$ means that the sequence $\{||X_m||\}$ is stochastically bounded, and \linebreak $X_m = o_p(1)$ means that $||X_m||\prob 0.$ $\mathrm{I}_p$ denotes the identity matrix of size $p.$
	
	\section{Observation model and consistency of the estimator}\label{s:2}
	\subsection{The EW-TLS promblem}
	We deal with the model $AX \approx B.$ Here $A\in\mathbb{R}^{m\times n}$ and $B\in\mathbb{R}^{m\times d}$ are matrices of observations, and the matrix $X\in\mathbb{R}^{n\times d}$ is to be estimated. Assume that
	\begin{equation}\label{eq:2.1}
	A = A_0 + \tilde{A}, \ \ B = B_0 + \tilde{B},
	\end{equation}
	and that there exists $X_0\in\mathbb{R}^{n \times d}$ such that
	\begin{equation}\label{eq:2.2}
	A_0 X_0 = B_0.
	\end{equation}
	Here $A_0$ is nonrandom true input matrix, $B_0$ is a true output matrix, and $\tilde{A},$ $\tilde{B}$ are error matrices. $X_0$ is the true value of the matrix parameter.
	
	It is useful to rewrite the model \eqref{eq:2.1} and \eqref{eq:2.2} as a classical errors-in-variables (EIV) model \cite{car}. Denote $a_i^{\top},$ $a_{0i}^{\top},$ $\tilde{a}_i^{\top},$ $b_i^{\top},$ $b_{0i}^{\top},$ $\tilde{b}_i^{\top},$ $i = 1,\dots, m,$ the rows of $A,$ $A_0,$ $\tilde{A},$ $B,$ $B_0$ and $\tilde{B},$ respectively. Then the model above is equivalent to the EIV model
	\begin{equation}\label{eq:2.3}
	a_i = a_{0i} + \tilde{a}_i, \ \ b_i = b_{0i} + \tilde{b}_i, b_{0i} = X_0^{\top}a_{0i}, \ \ i = 1,\dots,m.
	\end{equation}
	Vectors $a_{0i}$ are nonrandom and unknown, and vectors $\tilde{a}_i,$ $\tilde{b}_i$ are random errors. Based on observations $a_i,$ $b_i,$ $i = 1,\dots,m,$ one has to estimate $X_0.$
	
	Rewrite the model \eqref{eq:2.1} and \eqref{eq:2.2} in an implicit way. Introduce matrices
	\begin{equation}\label{eq:2.4}
	C = [A, B], \ \ C_0 = [A_0, B_0], \ \ \tilde{C} = [\tilde{A}, \tilde{B}], \ \ 
	Z_0 = 
	\begin{bmatrix}
	X_0 \\
	-\mathrm{I}_d
	\end{bmatrix}.
	\end{equation}
	Then \eqref{eq:2.1}, \eqref{eq:2.2} is equivalent to the next relations:
	\begin{equation*}\label{eq:2.5}
	C = C_0 + \tilde{C}, \ \ C_0 Z_0 = 0.
	\end{equation*}

	Let $\tilde{C} = (\tilde{c}_{ij}, i = 1,\dots,m, \ \ j = 1,\dots, n+d).$ Following \cite{kuhu} we state global assumptions of the paper, conditions \ref{i} to \ref{iv}.
	\renewcommand{\theenumi}{(\roman{enumi})}
	\renewcommand{\labelenumi}{(\roman{enumi}).}
	\begin{enumerate}
		\item\label{i}
		Vectors $\tilde{c}_i := (\tilde{c}_{i1},\dots,\tilde{c}_{i,n+d})^{\top},$ $i = 1,2,\dots,$ are independent with zero mean and finite second moments.
	\end{enumerate}
	
	Let $\sigma_{ij}^2 = \M \tilde{c}_{ij}^2,$ $i = 1,2,\dots,$ $j=1,\dots, n+d.$ We allow that some of $\sigma_{ij}^2$ are vanishing.
	\begin{enumerate}\addtocounter{enumi}{1}
		\item\label{ii} 
		For a fixed $J\subset \{1,2,\dots,n+d\},$ every $j\notin J$ and every $i = 1,2,\dots $ satisfy $\sigma_{ij}^2 = 0.$ Moreover
		$$
		\cov(\tilde{c}_{ij}, j\in J) = \sigma^2\Sigma_i, \ \ i = 1,2,\dots,
		$$
		with unknown positive factor of proportionality $\sigma^2$ and known matrices $\Sigma_i.$
		\item\label{iii} 
		There exists $\varkappa  > 0$ such that for every $i = 1,2,\dots,$ it holds $\lambda_{min}(\Sigma_i)\geq\varkappa^2.$ 
	\end{enumerate}	
	
	For the matrix $Z_0 = (z_{0,jk})$ given in \eqref{eq:2.4} and the set $J$ from condition \ref{ii}, denote
	\begin{equation*}\label{eq:2.6}
	Z_{0J} = (z_{0,jk}, j\in J, k =1,\dots, d).
	\end{equation*}
	\begin{enumerate}\addtocounter{enumi}{3}
		\item\label{iv}
		\begin{equation*}\label{eq:2.7}
		\mathrm{rank}(Z_{0J}) = d.
		\end{equation*} 
	\end{enumerate}	
	
	The EW-TLS problem consists in finding the value $\hat{X}$ of the unknown matrix $X$ and values of disturbances $\Delta\hat{A},$ $\Delta\hat{B}$ minimizing the weighted
	sum of squared corrections:
	\begin{equation}\label{eq:2.8}
	\min_{(X\in\mathbb{R}^{n\times d}, \Delta A, \Delta B)}\sum_{i=1}^{m}||\Sigma_i^{-1/2}\Delta c_i^J||^2
	\end{equation}
	subject to constrains
	\begin{equation*}\label{eq:2.9}
	(A - \Delta A)X = B - \Delta B, \ \ \Delta c_i^J = 0, \ \  i=1,\dots,m, \ \ j\notin J.
	\end{equation*}
	Here $C = [A,B] = (c_{ij}),$ $\Delta C = [\Delta A,\Delta B] = (\Delta c_{ij})$ and the column vectors
	\begin{equation*}\label{eq:2.10}
	\Delta c_i^J := (\Delta c_{ij}, j\in J)\in\mathbb{R}^{|J|}.
	\end{equation*}
	
	\subsection{EW-TLS estimator and its consistency}
	For a random realization, it can happen that the problem \eqref{eq:2.8}  has no solution. Assume conditions \ref{i} -- \ref{iv}.
	\begin{dfn}\label{def:1}
		The EW-TLS estimator $\hat{X} = \hat{X}_{EW-TLS}$ of $X_0$ in the model \eqref{eq:2.1}, \eqref{eq:2.2} is a Borel measurable mapping of the data matrix $C$ into $\mathbb{R}^{n\times d}\cup\{\infty \},$ which solves the problem \eqref{eq:2.8} under the additional constraint
		\begin{equation}\label{eq:2.11}
		rank(Z_J) = d
		\end{equation}
		$\left(\text{here } Z = 
		\begin{bmatrix}
		X \\
		-\mathrm{I}_d
		\end{bmatrix} = 
		(z_{jk}),\: Z_J:= (z_{jk}, j\in J, k = 1,\dots,d)\right),$ if there exists a solution, and $\hat{X} = \infty$ otherwise.
	\end{dfn}
	
	The EW-TLS estimator always exists due to \cite{pf}. We need more conditions to provide the consistency of $\hat{X}.$
	\begin{enumerate}\addtocounter{enumi}{4}
		\item\label{v}
		There exists $r\geq 2$ with $r> d\left(|J| - \dfrac{d+1}{2}\right)$ such that
		$$
		\sup_{(i\geq 1, j\in J)}\M|\tilde{c}_{ij}|^{2r} < \infty.
		$$
		\item\label{vi} 
		$\dfrac{\lambda_{\mathrm{min}}(A_0^{\top}A_0)}{\sqrt{m}} \to\infty,$ as $m\to\infty.$
		\item\label{vii} 
		$\dfrac{\lambda_{\min}^2(A_0^{\top}A_0)}{\lambda_{\max}(A_0^{\top}A_0)} \to\infty,$ as $m\to\infty.$
	\end{enumerate}
	
	The next result on weak consistency is stated in Theorem 2, \cite{kuhu}.
	
	\begin{thm}\label{th:2}
		Assume conditions \ref{i} to \ref{vii}. Then the EW-TLS estimator $\hat{X}$ is finite with probability tending to one, and $\hat{X}$ tends to $X_0$ in probability, as $m\to\infty.$
	\end{thm}
	
	Notice that under a bit stronger assumptions on eigenvalues of $A_0^{\top}A_0,$ the estimator $\hat{X}$ is strongly consistent, see Theorem 3, \cite{kuhu}.
	\subsection{The  estimating function}
	Remember that error vectors $\tilde{c}_i$ enter condition \ref{i} and the matrix $Z = Z(X)$ is introduced in Definition \ref{def:1}. Let
	\begin{equation*}\label{eq:2.12}
	S_i := \frac{1}{\sigma^2}\cov(\tilde{c}_i), i=1,2,\dots
	\end{equation*}
	
	Denote also
	\begin{equation}\label{eq:2.13}
	q(c,S;X) = c^{\top}Z(Z^{\top}SZ)^{-1}Z^{\top}c,
	\end{equation}
	where 
	$c = \begin{bmatrix}
	a\\
	b
	\end{bmatrix}
	\in\mathbb{R}^{(n+d)\times 1},$ $S\in \mathbb{R}^{(n+d)\times (n+d)},$ and  
	\begin{equation}\label{eq:2.14}
	Q(X) = \sum_{i=1}^m q(c_i, S_i; X), \ \ X\in\mathbb{R}^{n\times d}, \ \ \mathrm{rank}(Z_J) = d.
	\end{equation}
	
	Notice that due to \ref{iv} $|J|\geq d,$ and under constraint \eqref{eq:2.11} $Z_J$ is of full rank. Then, under conditions \ref{i} -- \ref{iii} the matrix $Z^{\top}S_i Z$ is nonsingular, $i = 1,2,\dots$
	
	The EW-TLS estimator is known to minimize the objective function \eqref{eq:2.13}, see Theorem~1, \cite{kuhu}.
	\begin{lem}\label{l:3}
		Assume conditions \ref{i} to \ref{iv}. The EW-TLS estimator $\hat{X}$ is finite if, and only if, there exists an unconditional minimum of the function \eqref{eq:2.14}, and then $\hat{X}$ is a minimum point of this function.
	\end{lem}
	
	Introduce an estimating function related to the loss function \eqref{eq:2.13}:
	\begin{gather}
	s(a, b, S; X) = \tilde{s}\cdot(Z^{\top} S Z)^{-1}, \label{eq:2.15} \\
	\tilde{s}=\tilde{s}(a, b, S; X) := ac^{\top}Z - [S_a, S_{ab}]Z(Z^{\top}SZ)^{-1}Z^{\top}cc^{\top}Z. \label{eq:2.15A}
	\end{gather}
	Here
	\begin{equation}\label{eq:2.15a}
	c = 
	\begin{bmatrix}
	a\\
	b
	\end{bmatrix},
	\ \ a\in\mathbb{R}^{n\times 1}; \ \  
	S = 
	\begin{bmatrix}
	S_a & S_{ab}\\
	S_{ba} & S_b
	\end{bmatrix},
	\ \ S_a\in\mathbb{R}^{n\times n}.
	\end{equation}
	\begin{cor}\label{cor:4}
		Assume conditions \ref{i} -- \ref{vii}. Then the next two statements hold true.
		\renewcommand{\theenumi}{(\alph{enumi})}
		\renewcommand{\labelenumi}{(\alph{enumi})}
		\begin{enumerate}
			\item\label{4.a}
			With probability tending to one $\hat{X}$ is a solution to the equation
			$$
			\sum_{i=1}^m s(a_i,b_i, S_i;X) = 0, \ \ X\in\mathbb{R}^{n\times d}, \mathrm{rank}(Z_J) = d.
			$$
			\item\label{4.b} 
			The function \eqref{eq:2.15} is an unbiased estimating function, i.e., for each $i\geq 1,$ $\M_{X_0}s(a_i,b_i,S_i;X_0) = 0.$
		\end{enumerate}
	\end{cor}
	
	For fixed $a,$ $b,$ $S,$ the function \eqref{eq:2.15} maps $X$ into $\mathbb{R}^{n\times d}.$ The derivative $s^{\prime}_X$ is a linear operator in this space.
	\begin{lem}\label{l:5}
		Under conditions \ref{i} -- \ref{vii}, for each $H\in\mathbb{R}^{n\times d}$ and $i\geq 1$ it holds
		\begin{equation}\label{eq:2.16}
		\M_{X_0}[s^{\prime}_X(a_i,b_i,S_i;X_0)\cdot H] = a_{0i}a_{0i}^{\top} H (Z_0 S_i Z_0)^{-1}.
		\end{equation}
	\end{lem}
	
	\section{Asymptotic normality of the estimator}\label{s:3}
	Introduce further assumptions.
	\begin{enumerate}\addtocounter{enumi}{7}
		\item\label{viii}
		For some $\delta>0,$ $\displaystyle\sup_{(i\geq 1, j\in J)}\M|\tilde{c}_{ij}|^{4 + 2\delta} < \infty.$
		\item\label{ix}
		For $\delta$ from the condition \ref{viii},
		$$
		\frac{1}{m^{1 + \delta/2}}\sum_{i=1}^m ||a_{0i}||^{2 + \delta} \to 0, \ \ \text{as } m\to\infty.
		$$
		\item\label{x}
		$\dfrac{1}{m} A_0^{\top} A_0 \to V_A,$ as $m\to\infty,$ where $V_A$ is a nonsingular matrix.
	\end{enumerate}
	
	Notice that condition \ref{x} implies assumptions \ref{vi}, \ref{vii}.
	
	\begin{enumerate}\addtocounter{enumi}{10}
		\item\label{xi}
		For matrices from condition the \ref{ii}, $\Sigma_i \to\Sigma_{\infty},$ as $m\to\infty,$ where $\Sigma_{\infty}$ is certain matrix.
	\end{enumerate}
	
	Notice that conditions \ref{xi}, \ref{iii} imply that $\Sigma_{\infty}$ is nonsingular.
	\begin{enumerate}\addtocounter{enumi}{11}
		\item\label{xii}
		If $p, q, r\in J$ (they are not necessarily distinct) and $i\geq 1,$ then 
		$$
		\M\tilde{c}_{ip}\tilde{c}_{iq}\tilde{c}_{ir} = 0.
		$$
		\item\label{xiii}
		If $p, q, r, u\in J$ (they are not necessarily distinct), then $\displaystyle\frac{1}{m}\sum_{i=1}^m \M\tilde{c}_{ip}\tilde{c}_{iq}\tilde{c}_{ir}\tilde{c}_{in}$ converges to a finite limit $\mu_4(p,q,r,u),$ as $m$ tends to infinity.
	\end{enumerate}
	
	Introduce a random element in the space of couples of matrices:
	\begin{equation}\label{eq:3.1}
	W_i = (a_{0i}\hat{c}_i^{\top}, \tilde{c}_i\tilde{c}_i^{\top} - \sigma^2 S_i).
	\end{equation}
	
	Hereafter $\convdistr$ stands for the convergence in distribution.
	\begin{lem}\label{l:6}
		Assume conditions \ref{i}, \ref{ii} and \ref{viii} -- \ref{xiii}. Then 
		\begin{equation}\label{eq:3.2}
		\frac{1}{\sqrt{m}}\sum_{i=1}^m W_i \convdistr \Gamma = (\Gamma_1,\Gamma_2), \ \ \text{as } m\to\infty,
		\end{equation}
		where $\Gamma$ is a Gaussian centered random element with independent matrix components $\Gamma_1$ and $\Gamma_2.$	
	\end{lem}
	
	Now, we state the asymptotic normality of the EW-TLS estimator.
	\begin{thm}\label{th:7} 
		Assume conditions \ref{i} -- \ref{v} and \ref{viii} -- \ref{xiii}. Then
		\begin{gather}
		\sqrt{m}(\hat{X} - X_0) \convdistr V_A^{-1}\Gamma(X_0), \text{as } m\to\infty, \label{eq:3.3} \\
		\Gamma(X) := \Gamma_1 Z + P_a\Gamma_2 Z - [S_a^{\infty}, S_{ab}^{\infty}]Z(Z^{\top}S_{\infty}Z)^{-1}(Z^{\top}\Gamma_2 Z), \label{eq:3.4} 
		\end{gather}
		where $V_A$ enters condition \ref{x}, $P_a$ is the projector with 
		$P_a
		\begin{bmatrix}
		a\\
		b
		\end{bmatrix}
		= a,$
		$\Gamma_1$ and $\Gamma_2$ enter relation \eqref{eq:3.2}, and
		\begin{equation}\label{eq:3.5}
		S_{\infty} = 
		\begin{bmatrix}
		S_a^{\infty} & S_{ab}^{\infty} \\
		S_{ba}^{\infty} & S_b^{\infty}
		\end{bmatrix}
		= \lim_{i\to\infty} S_i, \ \
		Z = 
		\begin{bmatrix}
		X\\
		-\mathrm{I}_d
		\end{bmatrix}.
		\end{equation}
		Moreover the limiting random matrix $X_{\infty} := V_A^{-1}\Gamma(X_0)$ has a nonsingular covariance structure, i.e., for each nonzero vector $u\in\mathbb{R}^{d\times 1},$ $\cov(X_{\infty}u)$ is a nonsingular matrix.
	\end{thm}
	
	\section{Construction of confidence region for a linear functional of $X_0$}\label{s:4}
	\subsection{Estimation of nuisance parameters}
	Theorem \ref{th:7} can be applied, e.g., to construct a confidence region for a linear functional of $X_0.$ For this purpose one has to estimate consistently a covariance structure of the limiting random matrix $V_A^{-1}\Gamma(X_0).$ Such a structure, besides of $X_0,$ depends on nuisance parameters. Some of them can be estimated consistently.
	
	Hereafter bar means average for rows $i = 1,\dots,m,$ e.g., 
	$$\overline{a b^{\top}} = m^{-1}\cdot\displaystyle\sum_{i=1}^m a_ib_i^{\top}, \ \ \bar{S} = m^{-1}\displaystyle\sum_{i=1}^m S_i.$$ 
	\begin{lem}\label{l:8}
		Assume conditions of Theorem \ref{th:7}. Define
		\begin{gather}
		\hat{Z} = 
		\begin{pmatrix}
		\hat{X}\\
		-\mathrm{I}_d
		\end{pmatrix}, \ \
		\hat{\sigma}^2 = \frac{1}{d}\mathrm{tr}\left[(\hat{Z}^{\top} \overline{c c^{\top}} \hat{Z}) (\hat{Z}^{\top} \bar{S} \hat{Z})^{-1}\right], \label{eq:4.1}	\\
		\hat{V}_A = \overline{a a^{\top}} - \hat{\sigma}^2\bar{S}. \label{eq:4.2}
		\end{gather}
		Then, as $m\to\infty,$
		\begin{equation*}\label{eq:4.3}
		\hat{\sigma}^2\prob \sigma^2, \ \ \hat{V}_A\prob V_A.
		\end{equation*}
	\end{lem}
	
	\subsection{Estimation of the asymptotic covariance structure of $X_0$}
	Let $u\in\mathbb{R}^{d\times 1},$ $u\neq 0.$ Theorem \ref{th:7} implies the convergence
	\begin{equation}\label{eq:4.4}
	\sqrt{m}(\hat{X}u - X_0 u)\convdistr N(0, S_u), \ \ \text{as } m\to\infty,
	\end{equation}
	with nonsingular matrix $S_u = \cov(V_A^{-1}\Gamma(X_0)u).$
	
	We start with the case of normal errors $\tilde{c}_i,$ $i = 1,2,\dots$ Then condition \ref{xii} holds true, and Theorem \ref{th:7} is applicable. The asymptotic covariance matrix $S_u$ is a continuous function $S_u = S_u(X_0, V_A, \sigma^2,S_{\infty})$ of unknown parameters (here the limiting covariance matrix $S_{\infty}$ could be unknown, though for a given $m,$ matrices $S_1,\dots, S_m$ are assumed known). Due to Theorem \ref{th:2} and Lemma \ref{l:8} the matrix 
	\begin{equation}\label{eq:4.5}
	\hat{S}_u:= S_u(\hat{X}, \hat{V}_A, \hat{\sigma}^2, \bar{S})
	\end{equation} 
	is a consistent estimator of $S_u.$ 
	
	Now, we do not assume the normality of the errors. Then the exact formula for $S_u$ does not allow to estimate it consistently, because the formula involves higher moments of errors which are difficult to estimate consistently. Instead, we use Corollary \ref{cor:4} to construct the so-called sandwich estimator \cite{car} for $S_u.$ Denote 
	\begin{equation}\label{eq:4.5a}
	\hat{s}_i = \tilde{s}(a_i, b_i, S_i; \hat{X}), \ \ i = 1,\dots,m,
	\end{equation}
	with $\tilde{s}$ introduced in \eqref{eq:2.15A} 
	\begin{lem}\label{l:9}
		Assume conditions of Theorem \ref{th:7}. For $u\in\mathbb{R}^{d\times 1},$ $u\neq 0,$ define 
		\begin{equation}\label{eq:4.6}
		\hat{S}_u = \hat{V}_A^{-1}\cdot\frac{1}{m}\sum_{i=1}^m \hat{s}_i u u^{\top}\hat{s}_i^{\top},
		\end{equation}
		with $\hat{V}_A$ given in \eqref{eq:4.2}, \eqref{eq:4.1}. Then $\hat{S}_u\prob S_u,$ as $m\to\infty.$  
	\end{lem}

	\textit{Remark}. In the case of normal errors, the estimator \eqref{eq:4.5} is asymptotically more efficient than the estimator \eqref{eq:4.6}, cf. the discussion in \cite{car}, p.~369.

	Given a consistent estimator $\hat{S}_u$ of $S_u,$ we have from \eqref{eq:4.4} that 
	\begin{equation}\label{eq:4.7}
	\sqrt{m}(\hat{S}_u)^{-1/2}(\hat{X}u - X_0 u)\convdistr N(0,\mathrm{I}_n), \ \ \text{as } m\to\infty.
	\end{equation}
	Based on \eqref{eq:4.7}, one can construct in a standard way an asymptotic confidence ellipsoid for $X_0 u.$ Similarly a confidence ellipsoid can be constructed for any finite set of linear combinations of $X_0$ entries.
	
	\section{Conclusion}\label{s:5}
	We proved the asymptotic normality of the EW-TLS estimator in a multivariate errors-in-variables model $AX\approx B$ with heteroscedastic errors.  We assumed the convergence \ref{xi} of the second error moments, vanishing third moments \ref{xiii}, and the convergence of averaged fourth moments \ref{xiii}. The condition \ref{xii} ensured that the asymptotic covariance structure of $\hat{X}$ is nonsingular. This condition holds  true in two cases: (a) all  the error vectors $\tilde{c}_i$ are symmetrically distributed, or (b) for each $i,$ random variables $\tilde{c}_{ip},$ $p\in J,$ are independent and have vanishing coefficient of asymmetry. 
	
	The obtained asymptotic normality result made it possible to construct a confidence ellipsoid for a linear functional of $X_0.$ Another plausible application is goodness-of-fit test in the model $AX \approx B$ with heteroscedastic errors (see \cite{kutsa} for such a test in the model with homoscedastic errors). 
	
	The author is grateful to Prof.~A.~Kukush for the problem statement and fruitful discussions. 
	\section*{Appendix}
	\renewcommand{\theequation}{A.\arabic{equation}}
	\subsection*{Proof of Corollary \ref{cor:4}}
	(a) The space  $\mathbb{R}^{n\times d}$ is endowed with natural inner product \linebreak $<A,B> = \mathrm{tr}(A B^{\top}).$ The matrix derivative $q_X^{\prime}$ of the functional \eqref{eq:2.13} is a linear functional on $\mathbb{R}^{n\times d},$ and based on the inner product, this functional can be identified with certain matrix from $\mathbb{R}^{n\times d}.$
	
	Remember that $Z = Z(X)$ is introduced in Definition \ref{def:1}. Using the rules of matrix calculus \cite{cart}, we have for $H\in\mathbb{R}^{n\times d}:$
	\begin{gather*}
	<q_X^{\prime}, H> = c^{\top}
	\begin{bmatrix}
	H\\
	0
	\end{bmatrix}
	(Z^{\top}S Z)^{-1} Z^{\top}c + c^{\top}(Z^{\top}S Z)^{-1}\cdot[H^{\top}, 0]c - \\
	- c^{\top}Z(Z^{\top} S Z)^{-1} \left([H^{\top}, 0]S Z + Z^{\top} S 
	\begin{bmatrix}
	H\\
	0
	\end{bmatrix}\right)
	(Z^{\top} S Z)^{-1} Z^{\top} c.
	\end{gather*}
	
	Remember relations \eqref{eq:2.15a}. Collecting similar terms, we obtain: 
	\begin{gather*}
	\frac{1}{2}<q_x^{\prime}, H> = a^{\top} H (Z^{\top} S Z)^{-1} Z^{\top}c - \\
	-c^{\top} Z (Z^{\top} S Z)^{-1}Z^{\top}
	\begin{bmatrix}
	S_a\\
	S_{ba}
	\end{bmatrix}
	H (Z^{\top} S Z)^{-1}Z^{\top} c,
	\end{gather*}
	and
	\begin{gather*}
	\frac{1}{2}<q_x^{\prime}, H> = \mathrm{tr}[ac^{\top} Z (Z^{\top} S Z)^{-1} H^{\top}] - \\
	\mathrm{tr}\left[ [S_a, S_{ab}]Z (Z^{\top} S Z)^{-1} Z^{\top} c c^{\top} Z (Z^{\top} S Z)^{-1} H^{\top}  \right]. 
	\end{gather*}
	
	Using the inner product in $\mathbb{R}^{n\times d}$ we obtain
	$$
	\frac{1}{2}q_x^{\prime} = \tilde{s}(X)(Z^{\top} S Z)^{-1},
	$$
	with $\tilde{s}(X) = \tilde{s}(a,b,S; X)$ given in \eqref{eq:2.15A}. Now, Theorem \ref{th:2} and Lemma \ref{l:3} imply the statement of Corollary \ref{cor:4}\ref{4.a}.
	
	(b) We set
	\begin{equation}\label{eq:A.1}
	a = a_0 + \tilde{a}, \ \ b = b_0 + \tilde{b}, \ \ b_0 = X^{\top}a_0, \ \ c = c_0 + \tilde{c} =   
	\begin{bmatrix}
	a_0\\
	b_0
	\end{bmatrix} +
	\begin{bmatrix}
	\tilde{a}\\
	\tilde{b}
	\end{bmatrix},
	\end{equation}
	where $a_0$ is a nonrandom vector and like in \eqref{eq:2.3},
	$$
	\cov(\tilde{c}) = \sigma^2 S = \sigma^2
	\begin{bmatrix}
	S_a    & S_{ab} \\
	S_{ba} & S_b
	\end{bmatrix}, \ \
	\M \tilde{c} = 0.
	$$
	Then
	\begin{gather*}
	\M_X a c^{\top} Z = a_0 c_0^{\top}
	\begin{bmatrix}
	X\\
	-\mathrm{I}_d
	\end{bmatrix} +
	\M \tilde{a}\tilde{c}^{\top} Z = \sigma^2[S_a, S_{ab}] Z,\\
	\M_X c c^{\top} Z = c_0 c_0^{\top}
	\begin{bmatrix}
	X\\
	-\mathrm{I}_d
	\end{bmatrix} +
	\M \tilde{c}\tilde{c}^{\top} Z = \sigma^2 S Z.
	\end{gather*}
	Therefore, see \eqref{eq:2.15},
	$$
	\M_X \tilde{s}(a,b,S;X) = \sigma^2[S_a, S_{ab}] Z - \sigma^2 [S_a, S_{ab}] Z (Z^{\top} S Z)^{-1} (Z^{\top} S Z) = 0
	$$
	The statement \ref{4.b} of Corollary \ref{cor:4} is proven. 
	
	\subsection*{Proof of Lemma \ref{l:5}}
	The derivative $s_X^{\prime}$ of the function \eqref{eq:2.15} with respect to $X$ is a linear operator in $\mathbb{R}^{n\times d}.$ Denote $f = f(Z) = Z(Z^{\top} S Z)^{-1}.$ For $H\in\mathbb{R}^{n\times d},$ it holds:
	$$
	\tilde{s}_X^{\prime} H = a a^{\top} H - [S_a, S_{ab}](f_X^{\prime} H)(Z^{\top}c c^{\top} Z) - [S_a, S_{ab}] f\cdot\left([H^{\top}, 0] c c^{\top} Z + Z^{\top} c c^{\top}
	\begin{bmatrix}
	H\\
	0
	\end{bmatrix}
	\right).
	$$
	
	We set \eqref{eq:A.1},  use relations 
	$$
	\M a a^{\top} = a_0 a_0^{\top} + \sigma^2S_a, \ \ \M_X (c c^{\top}Z) = \sigma^2 S Z
	$$
	and get:
	\begin{equation}\label{eq:A.2}
	\begin{array}{c}
	\M_X(\tilde{s}^{\prime}_X H) = (a_0 a_0^{\top} + \sigma^2 S_a)H -  
	\sigma^2[S_a, S_{ab}](f^{\prime}_X H)(Z^{\top} S Z) - \\
	- \sigma^2[S_a, S_{ab}]f\cdot\left([H^{\top}, 0] S Z + Z^{\top} S
	\begin{bmatrix}
	H\\
	0
	\end{bmatrix}
	\right).
	\end{array}
	\end{equation}
	
	Next,
	\begin{equation}\label{eq:A.3}
	f^{\prime}_XH = 
	\begin{bmatrix}
	H\\
	0
	\end{bmatrix}
	(Z^{\top} S Z)^{-1} - Z(Z^{\top} S Z)^{-1}\left(  [H^{\top}, 0]S Z + Z^{\top}S
	\begin{bmatrix}
	H\\
	0
	\end{bmatrix}  \right)
	(Z^{\top} S Z)^{-1}.
	\end{equation}
	
	Combining \eqref{eq:A.2} and \eqref{eq:A.3} we see that on the right-hand side of \eqref{eq:A.2} summands containing $H^{\top}$ are cancelled out.  We get finally
	$$
	\M_X(\tilde{s}_X^{\prime} H) = a_0 a_0^{\top} H, 
	$$
	which implies the statement, because by Corollary \ref{cor:4}\ref{4.b} it holds $\M_X\tilde{s}(a, b, S; X) = 0.$
	
	\subsection*{Proof of Lemma \ref{l:6}}
	The proof is similar to the proof of Lemmas 6 and 7 from \cite{kutsa} and based on Lyapunov's Central Limit Theorem. We just notice that due to condition \ref{xii} the matrix components of $W_i,$ namely $a_{0i}\tilde{c}_i^{\top}$ and $\tilde{c}_i \tilde{c}_i^{\top} - \sigma^2 S_i,$ are uncorrelated, and this implies the independence of matrix components $\Gamma_1$ and $\Gamma_2$ in \eqref{eq:3.2}.
	
	\subsection*{Proof of Theorem \ref{th:7}}
	We follow the line of \cite{kutsa}, see there the proof of Theorem 8(a). By Corollary \ref{cor:4}\ref{4.a}, it holds with probability tending to 1:
	\begin{equation}\label{eq:A.3a}
	\sum_{i=1}^m s(a_i, b_i, S_i; \hat{X}) = 0.
	\end{equation}
	
	Denote
	$$
	\hat{\Delta} = \sqrt{m}(\hat{X} - X_0), \ \ y_m = \sum_{i=1}^m s(a_i, b_i, S_i; X_0), \ \ U_m = \sum_{i=1}^m s_X^{\prime}(a_i, b_i, S_i; X_0).
	$$
	Using Taylor's formula around $X_0$ (see \cite{cart}, Theorem 5.6.2), we obtain from \eqref{eq:A.3a} that 
	\begin{equation}\label{eq:A.4}
	\begin{array}{c}
	\left(\dfrac{1}{m} U_m\right)\hat{\Delta} = -\dfrac{1}{\sqrt{m}} y_m + \mathrm{rest}_1,\\
	||\mathrm{rest}_1||\leq ||\hat{\Delta}||\cdot ||\hat{X} - X_0||\cdot O_p(1).
	\end{array}
	\end{equation}
	Here $O_p(1)$ is a multiplier of the form
	\begin{equation}\label{eq:A.5}
	\frac{1}{m}\sum_{i=1}^m \sup_{(||X - X_0||\leq \varepsilon_0)} ||s_X^{\prime\prime}(a_i, b_i, S_i; X)||,
	\end{equation}
	with positive $\varepsilon_0$ chosen such that $\mathrm{rank}(Z_J) = d$, for all $X$ with $||X - X_0||\leq \varepsilon_0;$ the choice is possible due to condition \ref{iv}, and expression \eqref{eq:A.5} is indeed $O_p(1)$ (i.e., stochastically bounded), because $s_X^{\prime\prime}$ is quadratic in $c_i$ and the averaged second moments of $c_i$ are assumed bounded. Thus, the relation \eqref{eq:A.4} holds true due to the consistency of $\hat{X}$  stated in Theorem \ref{th:2}.
	
	We have $||\mathrm{rest}_1||\leq ||\hat{\Delta}||\cdot o_p(1).$ Now, by Lemma \ref{l:5} and condition \ref{x} and \ref{xi} it holds 
	$$
	\frac{1}{m} U_m H = V_A H(Z_0^{\top} S_{\infty} Z_0)^{-1} + o_p(1), \ \ H\in\mathbb{R}^{n\times d},
	$$
	and we derive from \eqref{eq:A.4} the relation 
	\begin{equation}\label{eq:A.6}
	V_A \hat{\Delta}(Z_0^{\top} S_{\infty} Z_0)^{-1} = -\frac{1}{\sqrt{m}} y_m + \mathrm{rest}_2, \ \ ||\mathrm{rest}_2|| \leq ||\hat{\Delta}||\cdot o_p(1). 
	\end{equation}
	The summands in $y_m$ have zero expectation by Corollary \ref{cor:4}\ref{4.b}. Remember that $c_{0i} Z_0 = 0$ and the projector $P_a$ is introduced in Theorem \ref{th:7}. Then, see \eqref{eq:2.15}, 
	\begin{gather*}
	\tilde{s}(a_i, b_i, S_i; X_0) = (a_{0i} + \tilde{a}_i)\tilde{c}_i^{\top} Z_0 - [S_{ai}, S_{bi}] Z_0 (Z_0^{\top} S_i Z_0)^{-1}(Z_0^{\top}\tilde{c}_i\tilde{c}_i^{\top} Z_0), \\
	\tilde{s}(a_i, b_i, S_i; X_0) = W_{i1} Z_0 + P_a W_{i2} Z_0 - [S_{ai}, S_{bi}]Z_0 (Z_0^{\top} S_i Z_0)^{-1}(Z_0^{\top}W_{i2} Z_0).
	\end{gather*}
	Here $W_{ij}$ are components of \eqref{eq:3.1}. By Lemma \ref{l:6} it holds, see \eqref{eq:3.4} and condition \ref{xi}:
	\begin{equation}\label{eq:A.7}
	\frac{1}{\sqrt{m}}y_m \convdistr \Gamma(X_0)(Z_0^{\top} S_{\infty} Z_0)^{-1}, \ \ \text{as } m\to\infty.
	\end{equation}
	
	Now, relations \eqref{eq:A.6}, \eqref{eq:A.7} and nonsingularity of $V_A$ imply $\hat{\Delta} = O_p(1)$ and by Slutsky's lemma 
	$$
	V_A\hat{\Delta}(Z_0^{\top} S_{\infty} Z_0)^{-1} \convdistr \Gamma(X_0)(Z_0^{\top} S_{\infty} Z_0)^{-1}, \ \ \text{as } m\to\infty.
	$$
	This implies the desired convergence \eqref{eq:3.3} -- \eqref{eq:3.5}.
	
	Let $u\in\mathbb{R}^{d\times 1},$ $u\neq 0.$ By Lemma \ref{l:6} the components $\Gamma_1$ and $\Gamma_2$ are independent. We have
	\begin{gather*}
	\cov(\Gamma(X_0) u)\geq \cov(\Gamma_1 Z_0 u) = \lim_{m\to\infty}\frac{1}{m}\sum_{i=1}^m \M(a_{0i} \tilde{c}_i^{\top} Z_0 u u^{\top} Z_0^{\top}\tilde{c}_i a_{0i}^{\top}) = \\
	= \sigma^2 V_A[u^{\top}(Z_0^{\top} S_{\infty} Z_0) u],
	\end{gather*} 
	and the latter matrix is positive definite, because $V_A$ and $Z_0^{\top} S_{\infty} Z_0$ are positive definite under the conditions of Theorem 7. Therefore, $\cov(X_{\infty} u)$ is a positive definite matrix as well. 
	
	\subsection*{Proof of Lemma 8}
	We have 
	\begin{gather}
	\M c_i c_i^{\top} = c_{0i} c_{0i}^{\top} + \sigma^2 S_i, ~~
	Z_0^{\top}(\M c_i c_i^{\top}) Z_0 = \sigma^2 Z_0^{\top} S_i Z_0, \notag\\
	\sigma^2 = \frac{1}{d} \mathrm{tr} \left[(Z_0^{\top} \overline{c_i c_i^{\top}} Z_0) (Z_0^{\top} \bar{S} Z_0)^{-1}\right] +o_p(1). \label{eq:A.8}
	\end{gather}
	
	Relation \eqref{eq:A.8} and the convergence $\hat{Z} \prob Z_0$ imply the desired convergence $\hat{\sigma}^2 \prob \sigma^2,$ as $m\to\infty.$
	
	Next,
	\begin{gather*}
	\hat{V}_A = \M\overline{a a^{\top}} + o_p(1) - \hat{\sigma}^2\bar{S} = \overline{a_0 a_0^{\top}} + (\sigma^2 - \hat{\sigma}^2) \bar{S} + o_p(1),\\
	\hat{V}_A \prob \lim_{m\to\infty} \overline{a_0 a_0^{\top}} = V_A.
	\end{gather*}
	
	\subsection*{Proof of Lemma \ref{l:9}}
	Denote $\tilde{s}_i = \tilde{s}(a_i, b_i, S_i; X_0),$ $i = 1,2,\dots$ Then expansion \eqref{eq:A.6} implies that
	$$
	S_u = V_A^{-1}\cdot\frac{1}{m}\sum_{i=1}^m \tilde{s}_i u_i u_i^{\top} \tilde{s}_i^{\top} + o_p(1),
	$$
	and by Lemma \ref{l:8}
	\begin{gather*}
	S_u = \hat{V}_A^{-1}\cdot\frac{1}{m}\sum_{i=1}^m \tilde{s}_i u_i u_i^{\top} \tilde{s}_i^{\top} + o_p(1),\\
	\hat{S}_u - S_u = \hat{V}_A^{-1}\cdot\frac{1}{m}\sum_{i=1}^m (\hat{s}_i u_i u_i^{\top} \hat{s}_i^{\top} - \tilde{s}_i u_i u_i^{\top} \tilde{s}_i^{\top}) + o_p(1).
	\end{gather*}
	Then $\hat{S}_u - S_u \prob 0,$ as $m\to \infty,$ because $\hat{Z} \prob Z_0$ and $\overline{c c^{\top}} = O_p(1)$ (see formulas \eqref{eq:2.15}, \eqref{eq:2.15A} and \eqref{eq:4.5a}). Lemma \ref{l:9} is proven.

\end{document}